\documentclass[12pt,leqno]{amsart}
\newtheorem{theorem}{Theorem}[section]

\numberwithin{equation}{section}
\usepackage{amssymb,amsfonts,amsthm,setspace,indentfirst}
\usepackage[dvips]{graphics}
\usepackage{epsfig}
\usepackage[dvips]{graphics}
\graphicspath{ {images/} }
\usepackage{epsfig}
\usepackage{wrapfig}
\usepackage{caption}
\usepackage{subcaption}

\numberwithin{figure}{section}

\title[A note on the isoperimetric inequality in the plane]{ {\normalsize \bf{A note on the isoperimetric inequality in the plane}}}
\author[A. A. Shaikh and C. K. Mondal]{Absos Ali Shaikh and Chandan Kumar Mondal}

\begin{document}

\address{\noindent\newline Department of Mathematics,\newline University of
Burdwan, Golapbag,\newline Burdwan-713104,\newline West Bengal, India}
\email{aask2003@yahoo.co.in, aashaikh@math.buruniv.ac.in}
\email{chan.alge@gmail.com}

\begin{abstract}
It is well known that among all closed bounded curves in the plane with the given perimeter, the circle encloses the maximum area. There are many proofs in the literature. In this article we have given a new proof using some ideas of Demar\cite{DEM75}.
\end{abstract}
\date{\today}
\maketitle

\section{\textbf{Introduction}}
The isoperimetric inequality is a very well known geometric problem. The problem is to find a simple closed curve that encloses the maximum area than any other simple closed curve having the same length. From the ancient time it is known that the solution is the circle. There is a famous story about this problem. The story is dating back the problem to Vergil’s epic Latin poem The \emph{Aeneid}\cite{NOR07}, written in the period 29 BC to 19 BC. This
poem is about the foundation of the city of Carthage. The Greeks \cite{POR33}, were the first who considered this problem and tried to solved it.\\
\indent Many mathematicians from ancient time tried to solve this problem. It was believed that Zenodours \cite{COO40} in second century first gave a complete proof. Later Steiner, Euler, Weierstrass and many mathematicians proved it using different techniques. For complete history and some proofs see Blasj\"o \cite{BLA05}, Osserman \cite{OSS78}. \\
\indent Demar\cite{DEM75} in 1975 proved the isoperimetric inequality in the plane using some geometric arguments. He proved that if a closed continuous curve solves the isoperimetric inequality in the plane then any triangle formed by the intersection of two tangents at any two points and the chord joining them is isosceles. And he then showed that for every point on that curve at which a tangent can be drawn lies at a same distance from a fixed point inside the curve. In this article we only use the first idea of Demar\cite{DEM75} that all such triangles are isosceles and then we shall prove that only for the region having maximum area all such triangles with equal base are congruent.

\section{\textbf{Preliminaries}}
Our main target is to prove the following theorem.
\begin{theorem}[Isoperimetric inequality]\label{iso}
If $\gamma$ be a simple closed curve in the plane with length $L$ and bounds a region of area $A$ then
$$L^2\geq 4\pi A,$$
where the equality holds if and only if $\gamma$ is a circle.
\end{theorem}
It is well known that the solution is a convex region in the plane \cite{BLA05}, \cite{DEM75}. So only the convex regions will be focused here. We will also assume that the solution of the above inequality exists. For the existence of such problem, see book of Russel Benson \cite{BEN66}. Before proving the above theorem we will describe some notions.\\

\par Let $\gamma$ be a simple closed curve bounding a convex region $R$. Now standing at any point $A$ on the curve $\gamma$ and facing the region $R$, there is a right-hand and left-hand direction along $\gamma$ from $A$. A right-hand tangent at $A$ is a ray from $A$ tangent to $\gamma$ in the right-hand direction and similarly left-hand tangent can be defined \cite{KRA61}. Now the interior angle between these two tangents, i.e, the angle between right hand and left hand tangent at any point, measured in the sector containing $R$, is at most $\pi$ for a convex region \cite{DEM75}. \\
\begin{wrapfigure}{r}{0.5\textwidth}
\vspace{-20pt}
  \begin{center}
    \includegraphics[width=0.48\textwidth]{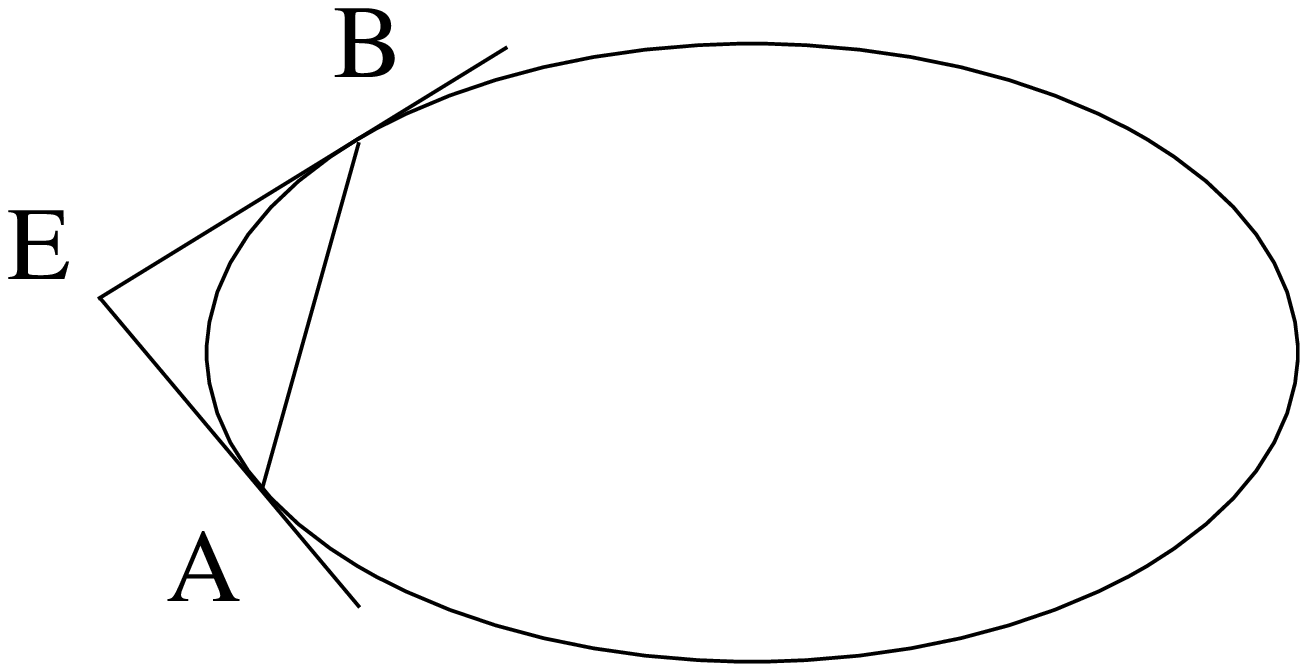}
  \end{center}
	\vspace{-40pt}
  \caption{}
  \label{fig1}
\end{wrapfigure}
\indent Consider two points $A$ and $B$ on $\gamma$ such that tangents at $A$ and $B$ meets each other at $E$. If the length $AB$ is $l$ then we call the triangle $\triangle ABE$, the tangent triangle on the curve $\widehat{AB}$ of base length $l$ and denoted it by $\triangle_l \widehat{AB}$ see Figure \ref{fig1}. It is not always possible that for any base length we get a tangent triangle, e.g, there is no tangent triangle of a circle with base length same as diameter. But for sufficiently small length it is possible to construct a tangent triangle at any two points on the curve by taking this length as base of the triangle. Now for general convex curve two tangent triangles with equal base length may not be congruent but in a circle all the tangent triangles with equal base are congruent.

\begin{proof}[Proof of theorem \ref{iso}]
We prove this by contradiction. Let $\gamma$ be a simple closed curve of length $L$, enclosing a convex region $R$. Let us assume that $R$ is not circular and has the property that it has maximum area than any other region with same boundary length $L$. Since the set of points where the convex curve is not differentiable is at most countable [\cite{HER55}, page 154] hence there is a dense set of points of $\gamma$ where tangents can be drawn. Demar \cite{DEM75} proved that any tangent triangle of $R$ is isosceles, see Figure \ref{fig2}. We will show that all the tangent triangles having same base length are congruent.\\

\indent Without loss of generality we can take two tangent triangles $\triangle_r\widehat{AB}$ and $\triangle_r\widehat{BC}$ which are not congruent. Since their base are equal hence $\angle EAB \neq \angle FBC$. Assume that $$\angle EAB >\angle FBC$$. Now the triangle $\triangle ABC$ is isosceles. So we get
$$\angle EAC>\angle FCA.$$
 Now reflect the whole area $ACFBE$ with respect to the perpendicular bisector of $AC$. Then we get a new curve $\gamma'$ by replacing $\widehat{AC}$ by its image $\widehat{A'C'}$ under the reflection. Hence we get a new region $R'$, see Figure \ref{fig3} and since reflection does not change area and perimeter hence both $R$ and $R'$ have same area and perimeter. Now the image of the triangles $\triangle ABE$, $\triangle CFB$ and $\triangle ABC$ under the reflection are $\triangle A'BE'$, $\triangle C'F'B$ and $\triangle BA'C'$ respectively.

		\begin{figure}[h!]
		
		\begin{center}
		\vspace{-0.0in}
		{\epsfxsize=10cm\epsfysize=8cm\epsfbox{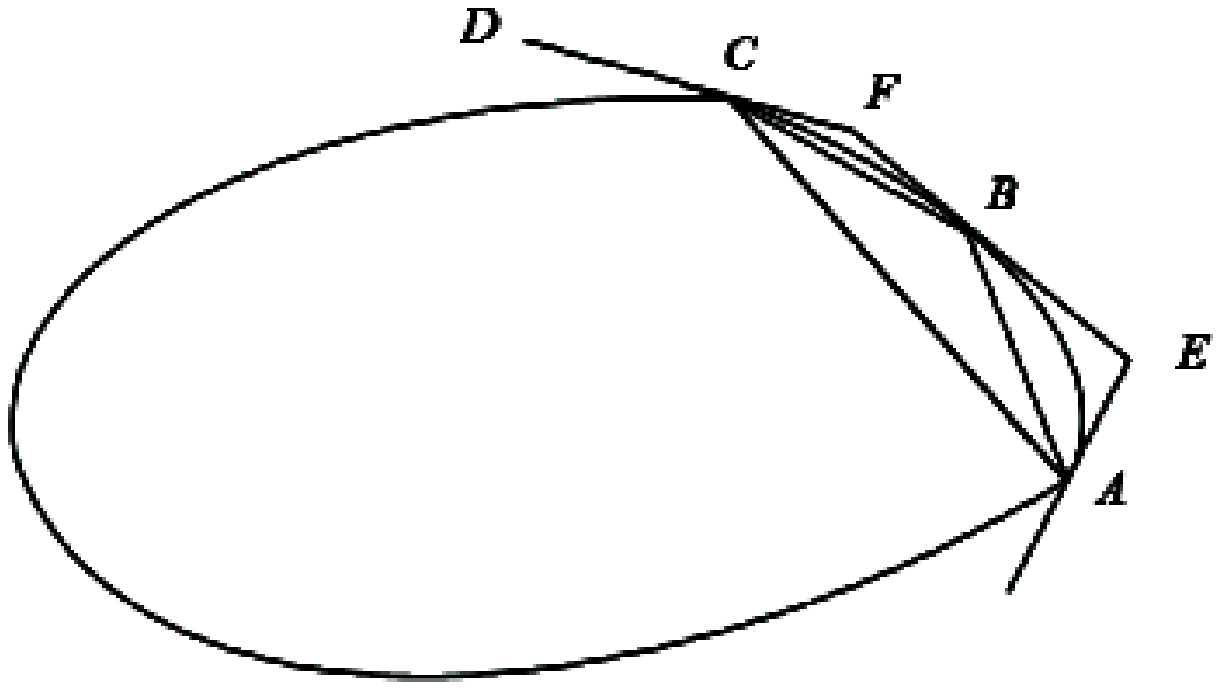}}
		\vspace{-0.5in}
		\caption{}
		\label{fig2} 
		\end{center}
		
		\end{figure}
		\vspace{-0.1in}
		%

  Hence we get $$\angle FCB=\angle F'BC'<\angle E'A'B.$$ 
This implies that $$\angle FCA=\angle F'C'A'<\angle E'A'C'=\angle EAC,$$
since $\triangle BA'C'$ is isosceles.\\
So combining above two inequalities we get
\begin{eqnarray*}
\pi&=&\angle FCA+\angle ACD,\\
& < & \angle E'A'C'+\angle C'A'D
\end{eqnarray*}

		\begin{figure}[h!]
		
		\begin{center}
		\vspace{-0.0in}
		{\epsfxsize=10cm\epsfysize=8cm\epsfbox{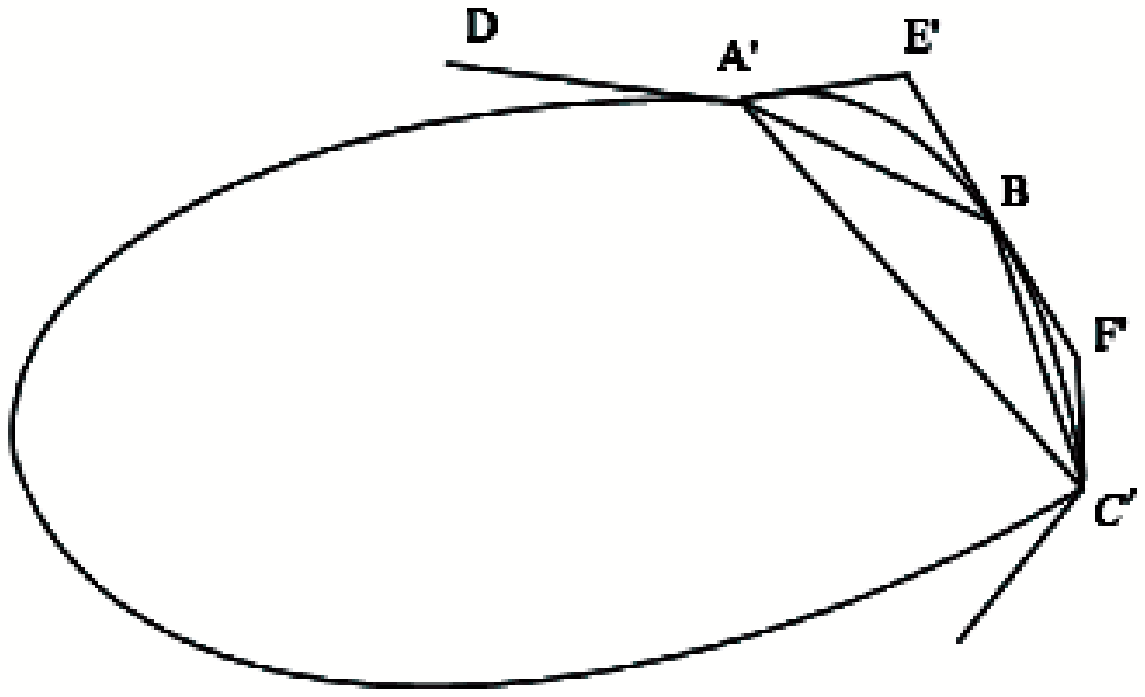}}
		\vspace{-0.5in}
		\caption{}
		\label{fig3} 
		\end{center}
		
		\end{figure}
		\vspace{-0.1in}
		%
But $A'E'$ and $A'D$ are right hand and left hand tangents of $\gamma'$, therefore $R'$ is not convex. It is a contradiction that $R$ is the region of maximum area among all regions with the same perimeter. 
Hence all the tangent triangles having equal base are congruent and since all the tangent triangles are isosceles, so $\gamma$ must be a circle.
\end{proof}

\end{document}